\documentclass{article}
\usepackage{acronym}

\newtheorem{theorem}{Theorem}
\newtheorem{acknowledgement}[theorem]{Acknowledgement}
\newtheorem{algorithm}[theorem]{Algorithm}

\newtheorem{example}[theorem]{Example}

\newenvironment{proof}[1][Proof]{\noindent\textbf{#1.} }{\ \rule{0.5em}{0.5em}}
\input{tcilatex}

\begin{document}

\title{The proof of Steinberg's three coloring conjecture}
\author{I. Cahit}
\date{email. \texttt{icahit@gmail.com}}
\maketitle

\begin{abstract}
The well-known Steinberg's conjecture asserts that any planar graph without $%
4$- and $5$-cycles is $3$ colorable. In this note we have given a short
algorithmic proof of this conjecture based on the spiral chains of planar
graphs proposed in the proof of the four color theorem by the author in 2004.
\end{abstract}

\section{ Introduction}

The problem which we will be dealt with in this note is known as "the three
color problem". Clearly the central problem is the four coloring of planar
graphs. Ruling out the almost trivial two coloring of bipartite graphs the
transition from three coloring to four coloring with respect to graphical
property e.g., planarity and non-planarity, cycle sizes, cycles through fix
number of vertices, connectivity (disjoint paths between two vertices) of
the graph is not very sharp. But there is an old attempt of Heawood for
three-color criterion of triangulated planar graphs which has been
strengthened in this paper as well as in [1].

Let us briefly list the previous results on the three coloring of planar
graphs. Apparently the first 3-color criterion for planar graphs was given
by Heawood in 1898 and is known as the \textquotedblleft three color
theorem"; simply saying that a finite triangulation is 3-colorable if and
only if it is even [1]. One may expect that after the proof of four color
theorem [2],[3] in 1976 and another earlier result due to Gr\"{o}tzsch
saying that \textquotedblleft all planar graphs without 3-cycles are three
colorable"\ [4] leads that all coloring problems related with the planar
graphs should be easily derived from these results . But Steinberg's three
colorability conjecture for planar graphs without 4 and 5 cycles tell us the
opposite. Erd\"{o}s [5] suggested an simpler question by asking whether
there exits a constant C such that the absence of cycles with size from 4 to
C in a planar graph guarantees its 3-colorability? In response to this
question first Abbott and Zhou [6] proved that such a C exists and C$\leq $%
11. Then Borodin improved it first to 10 and then to C$\leq $9 [7],[8].
Independently Sanders and Zhao have obtained the same upper bound [9]. Very
recently Borodin et. al. reduce the upper-bound to 7 by showing that [10]:

\begin{theorem}
Every planar graph without cycles of length 4 to 7 is 3-colorable.
\end{theorem}

In fact they showed a little stronger result:

\begin{theorem}
Every proper 3-coloring of the vertices of any face of size from 8 to 11 in
a connected graph G$_{7}$ can be extended to a proper 3-coloring of the
whole graph.
\end{theorem}

In Theorem 2, G$_{7}$ denote the class of planar graphs without cycles of
size from 4 to 7.

Most of these results on this conjecture are based on the discharging method
which was first used in the proof of the four-color theorem. Similarly let
us denote by G$_{6}$ , the class of planar graphs without cycles of size
from 4 to 6. On the other hand in [11] a new 3-color criterion has been
given for planar graphs:

\begin{theorem}
(3-Color Criteria). Let G be a (biconnected) plane graph. The following
three conditions are equivalent:
\end{theorem}

(i) G is 3-colorable.

(ii) There exists an even triangulation H$\supseteq $G.

(iii) G is edge-side colorable.

\bigskip

\section{Spiral chain coloring algorithm}

Our method is not rely on the previous results or necessary conditions of $3$
colorability of planar graphs. In fact we have given a coloring algorithm
which color vertices of any planar graph in $G_{6}.$ Let us assume that the $%
G_{6}$ has been embedded in the plane without edge crossing. Since the graph
G$_{6}$ has no cycle of length four and five this implies that it has also
no subgraph such as two triangles with a common edge.

Let $c_{1},c_{2},c_{3}$ be the three colors, say green, yellow, and red. Let
us define spiral-chains in $G_{6}$ which is exactly same as in [12],[13].
That we select any vertex on the outer-cycle of $G_{6}$ and scan outer
leftmost vertices in the clockwise (or counter clockwise) direction and
bypassing a vertex that is already scanned. If the last scanned vertex is
adjacent only previously scanned vertices then we select closest vertex to
the last scanned vertex and begin for a new spiral-chain till all vertices
scanned. This way we have obtained the set of vertex disjoint spiral chains $%
S_{1},S_{2},...,S_{k}$ of $G_{6}$. Note that if $\ k=1$ then $S_{1}$ is also
an hamilton path of $G_{6}$. Coloring of the vertices in $%
S_{1},S_{2},...,S_{k}$ are carried out in the reverse order of the spiral
chains and vertices i.e., $S_{k},S_{k-1},...,S_{1}.$ A very simple
spiral-chain coloring algorithm is given below which is similar the one
given for edge coloring of cubic planar graphs in [13] but here we do not
need Kempe switching.

\begin{algorithm}
Let $S_{1},S_{2},...,S_{k}~$\ be the spiral chains of G$_{6}$ , where spiral
chains are ordered from outer-edges towards inner edges.

Then color respectively the vertices of $S_{k},S_{k-1},...,S_{2},S_{1}$
using smaller indexed colors whenever possible, that is $c_{1}$ (green) is
in general the most frequent and $c_{3}$ (red) is the least frequent color
used in the algorithm. Note that while coloring the vertices of spiral chain 
$S_{i},k\geq i\geq 1$ if vertices $v_{j}$ and $v_{j+1}$ of an edge $%
(v_{j},v_{j+1})\in S_{i}$ is also belongs to a triangle and has been
colored, say with $c_{1}$ and $c_{2}$ then the other vertex $v_{k}$ of the
triangle will get the color $c_{3}$. If $v_{k}$ has already been colored,
say by $c_{2}$ then $v_{j}$ is colored by $c_{1}$ and $v_{j+1}$ is colored
by $c_{3}$.
\end{algorithm}

Cases will be studied in the proof of the following theorem.

\begin{theorem}
Algorithm 4 colors the vertices of G$_{6}$ with at most three colors.
\end{theorem}

\begin{proof}
By the ordering of spiral-chains that if $S_{i}\succ S_{j}$ then vertices of 
$S_{i}$ are colored before the vertices of $S_{j}$ and also if $v_{i}\succ
v_{j}$ and $v_{i},v_{j}\in S_{k}$ then the vertex $v_{i}$ is colored before $%
v_{j}$. Worst case configuration (a subgraph) in $G_{6}$ is a cycle of
length six with triangles attached to its edges is shown in Fig.1. The first
spiral chain $S_{1}$ passing through the three vertices forces all other
four spiral chains $S_{2},S_{3},S_{4},S_{5}$ ordered in the clockwise
direction. Now the spiral-coloring starts from $S_{5}$ and go to $S_{1}$%
without needing the forth color as follows:

$S_{5}:$ $%
c_{1},c_{2},...,S_{4}:c_{2},c_{3},...,S_{3}:c_{1},c_{3},...,S_{2}:c_{2},c_{3},..., 
$ and $S_{1}:c_{2},c_{1},c_{3},...$. Note that because of six triangles we
have to use the color "red" six times but this is not the case for a less
complicated configuration i.e., with less number of triangles, which can be
obtained from this one by only deleting \ vertices colored red.

Consider three disjoint triangles $%
T_{1}=(x_{1},y_{1},z_{1}),T_{2}=(x_{2},y_{2},z_{2})$ and $%
T_{3}=(x_{3},y_{3},z_{3})$ and three spiral chains $S_{i},S_{j},S_{k},i\leq
j\leq k$ respectively passing through the edges $x_{i}y_{i},i=1,2,3$ . Let $%
v $ be a vertex such that $v\in S_{p},(z_{i},v)\in E(G_{6}),i=1,2,3.$
Without loss of generality assume that $p>k$ that is $S_{k}$ is coming after 
$S_{p}$ and $x_{1}y_{1}\in S_{i},x_{2}y_{2}\in S_{j},x_{3}y_{3}\in S_{k}$.
We will show that vertices $z_{1},z_{2},z_{3}$ cannot have three different
colors but only $c_{3}$ by the spiral chain coloring. Assume contrary that
the vertices $x_{i},y_{i},z_{i}$ have the following coloring:

$x_{3}\rightarrow c_{1},y_{3}\rightarrow c_{2}\Longrightarrow
z_{3}\rightarrow c_{2}$

$x_{2}\rightarrow c_{2},y_{2}\rightarrow c_{3}\Longrightarrow
z_{2}\rightarrow c_{1}$

$x_{1}\rightarrow c_{1},y_{1}\rightarrow c_{2}\Longrightarrow
z_{1}\rightarrow c_{3}.$

But by the coloring rule of the algorithm (piority of the colors in the
assignment) $z_{1},z_{2},z_{3}$ are forced to be colored by $c_{3}.$ So we
can start coloring the first vertex $v$ of $S_{p}$ with $v\rightarrow c_{1}$%
. Lastly its easy to see that no any vertex of $G_{6}$ are left not colored.
\end{proof}

Based on the algorithm and proof of the above theorem we state.

\begin{theorem}
Planar graphs without $4-$ and $5-$cycles are three colorable.
\end{theorem}

\begin{example}
In Figure 1 we illustrate our spiral-chain three coloring for a graph with
19 vertices in G$_{6}$. (a) First we select arbitrarily a vertex from the
outer-cycle of the graph and starting from that vertex and traversing the
vertices in clockwise direction (always select not yet scanned
outer-leftmost vertex) in the form of a spiral chain. Hence the spiral-chain 
$S_{1}$ with vertex set $V(S_{1})=\{v_{1},v_{2},...,v_{15}\}$. \ Since the
last vertex $v_{15}$ is only adjacent to a vertex of $S_{1}$ we start new
spiral chain $S_{2}$ from the closest unscanned adjacent vertex (in the
graph vertex $v_{16}$) to $v_{15}$. Similarly we have obtained the second
spiral chain $S_{2}$ with the vertex set $V(S_{2})=\{v_{16},v_{17},v_{18}\}$%
. \ Since the other edges connected to $v_{18}$ have adjacent to already
scanned spiral chain $S_{1}$ vertices ($v_{10}$ and $v_{11}$) we start a new
spiral-chain $S_{3\text{ }}$closest unscanned vertex which is a isolated
vertex $V(S_{3})=\{v_{19}\}$. (b) Second part of the algorithm is the
three-coloring of the vertices of the spiral-chains \thinspace $%
S_{3},S_{2},S_{1}$ in reveres order. Colors $c_{1}$ (green),$c_{2}$ (yellow)
and $c_{3}$ (red) have been denoted by the numbers $1,2$ and $3$. Note that
in the spiral-chain coloring algorithm the color "green" has a pioroty over
"yellow" and "red" colors and the color yellow has a piority over the red
color. Color the vertex in $V(S_{3})=\{v_{19}\}$ with $v_{19}\rightarrow
c_{1}$ and since $(v_{19},v_{14},v_{13})$ is a triangle we also color $%
v_{14}\rightarrow c_{2},v_{13}\rightarrow c_{3}$. Color the vertices of $%
V(S_{2})=\{v_{18},v_{17},v_{16}\}$ as $v_{18}\rightarrow
c_{1},v_{17}\rightarrow c_{2},v_{16}\rightarrow c_{1}$. Since $%
(v_{18},v_{11},v_{10})$ is a triangle we also color $v_{12}\rightarrow
c_{2},v_{10}\rightarrow c_{3}.$ Similarly since $(v_{16},v_{8,}v_{7})$ is a
triangle we color $v_{8}\rightarrow c_{2}$ and $v_{7}\rightarrow c_{3}$.
Finally color the vertices of $V(S_{1})=\{v_{15},v_{14},...,v_{1}\}$ as $%
v_{15}\rightarrow c_{1}$, since $(v_{15},v_{7},v_{6})$ is a triangle color
also $v_{7}\rightarrow c_{3},v_{7}\rightarrow c_{2}.$Since $v_{2}$ and $v_{3}
$ have been already colored with $c_{2}$ and $c_{3}$ we continiue coloring $%
v_{12}\rightarrow c_{1}$ and since $v_{11}$ and $v_{10}$ have been colored
before with $c_{2}$ and $c_{3}$ we color $v_{9}\rightarrow c_{1}.$Since $%
v_{8},v_{7}$ and $v_{6}$ are colored before with $c_{2},c_{3}$ and $c_{2}$
we color $v_{6}\rightarrow c_{1},v_{3}\rightarrow c_{1}$ and $%
v_{1}\rightarrow c_{1}$.
\end{example}

\begin{acknowledgement}
\bigskip I like to thank to C.\ C. Heckman for his interest and views on
this problem.
\end{acknowledgement}

\bigskip

\bigskip

\newpage

\FRAME{ftbpFU}{195.125pt}{193.25pt}{0pt}{\Qcb{Worst case three coloring
around cycle of length three. Coloring is carried out from high index spiral
chains toward low index spiral chains.}}{}{Figure}{\special{language
"Scientific Word";type "GRAPHIC";maintain-aspect-ratio TRUE;display
"USEDEF";valid_file "T";width 195.125pt;height 193.25pt;depth
0pt;original-width 803.125pt;original-height 795.25pt;cropleft "0";croptop
"1";cropright "1";cropbottom "0";tempfilename
'J2NDV000.wmf';tempfile-properties "XNPR";}}\FRAME{ftbpFU}{300.8125pt}{229pt%
}{0pt}{\Qcb{Illustration of the spiral-chain three coloring.}}{}{Figure}{%
\special{language "Scientific Word";type "GRAPHIC";maintain-aspect-ratio
TRUE;display "USEDEF";valid_file "T";width 300.8125pt;height 229pt;depth
0pt;original-width 778.8125pt;original-height 591.9375pt;cropleft
"0";croptop "1";cropright "1";cropbottom "0";tempfilename
'J2P8SU00.wmf';tempfile-properties "XPR";}}

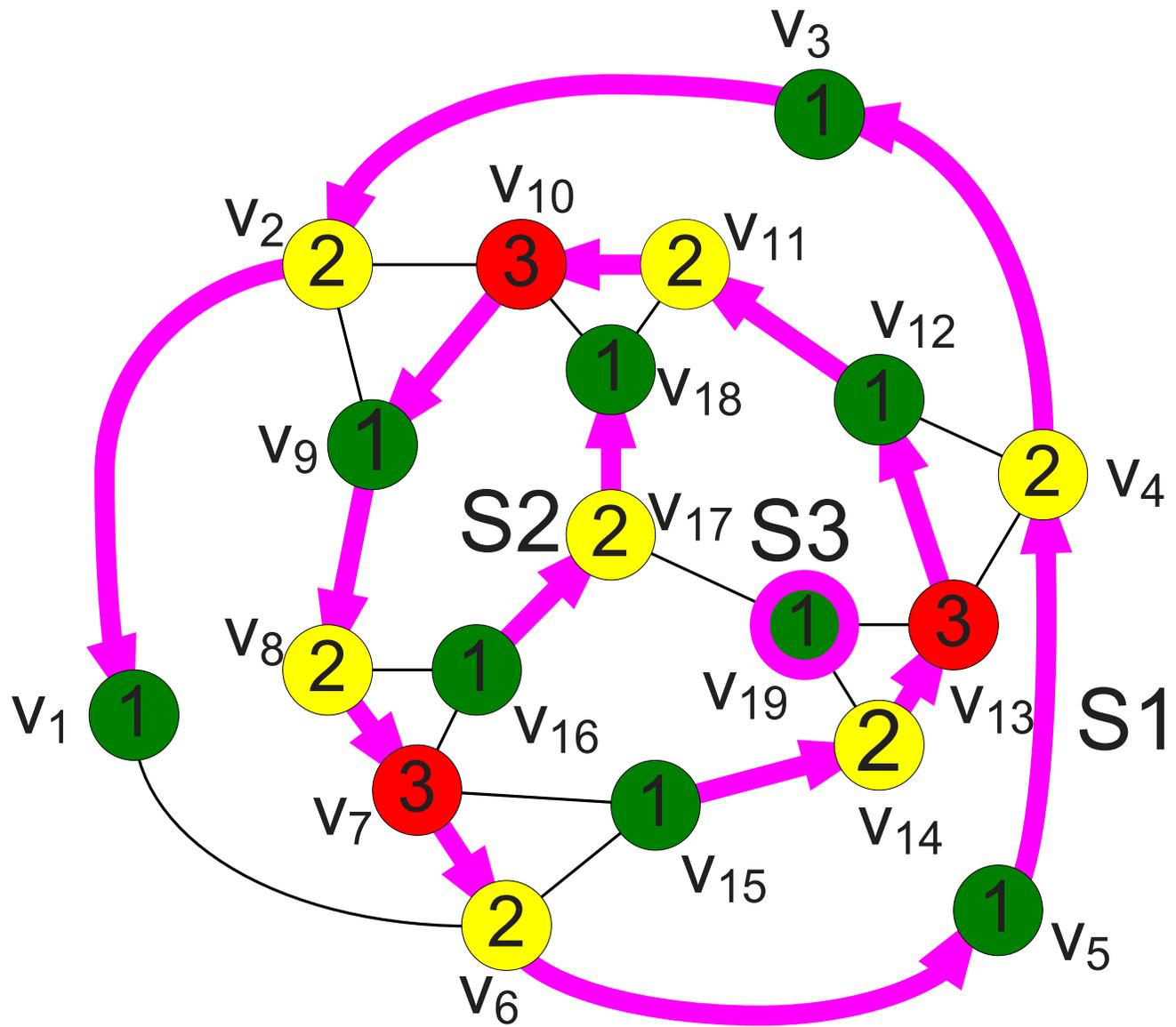
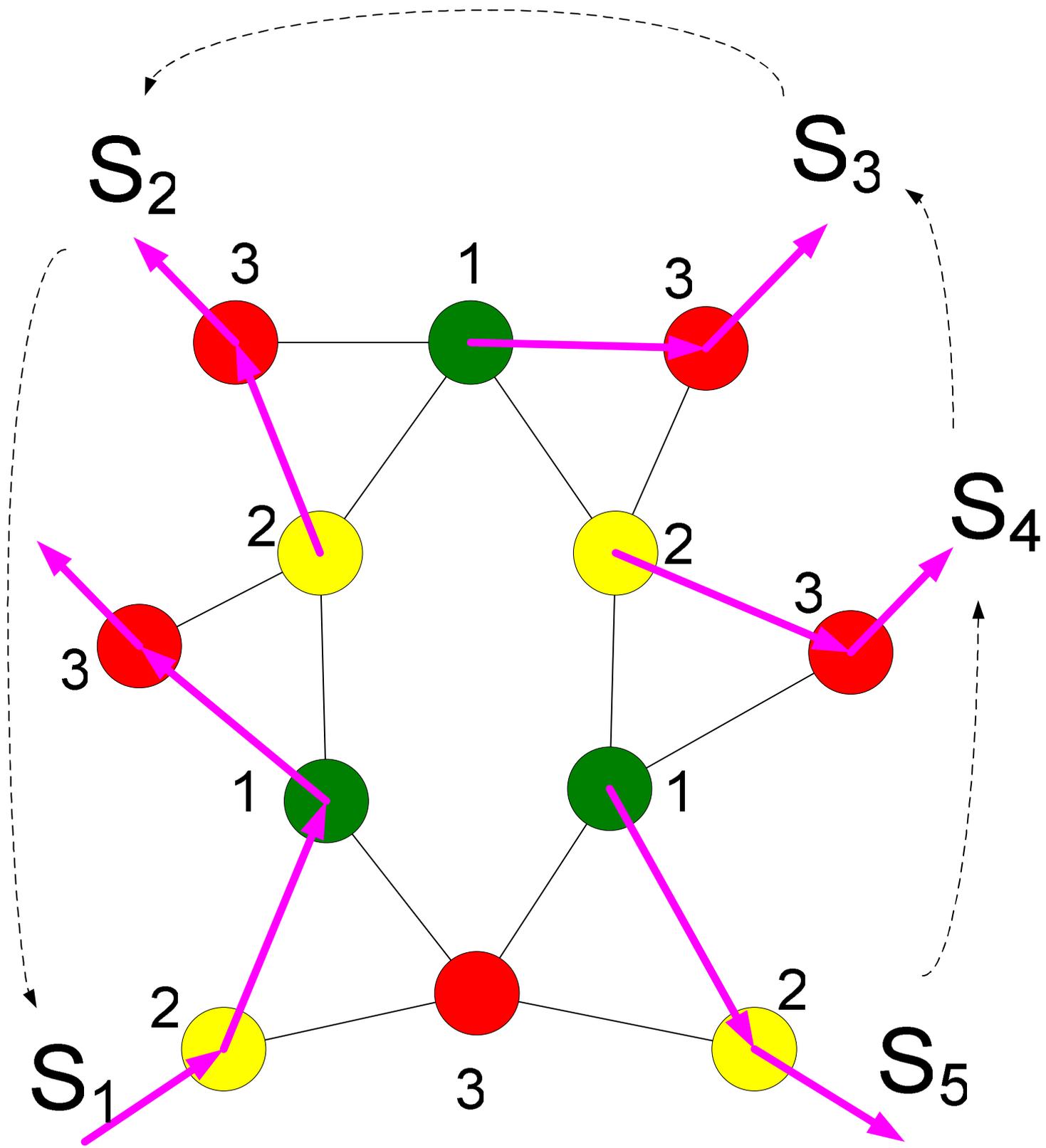
\end{document}